\documentclass[11pt]{amsart}
\usepackage{mathrsfs}
\usepackage{amssymb}

\usepackage{titletoc}
\usepackage{amscd}
\usepackage{amsmath, amssymb}
\usepackage{amsfonts}
\usepackage[colorlinks,linkcolor=blue,citecolor=blue, pdfstartview=FitH]{hyperref}

\numberwithin{equation}{section}

\theoremstyle{plain}
\newtheorem{thm}{Theorem}[section]
\newtheorem{theorem}[thm]{Theorem}
\newtheorem{lemma}[thm]{Lemma}
\newtheorem{corollary}[thm]{Corollary}
\newtheorem{proposition}[thm]{Proposition}

\theoremstyle{definition}

\newtheorem{remark}[thm]{Remark}

\newtheorem{definition}[thm]{Definition}

\newtheorem{example}[thm]{Example}

\newtheorem{defn-thm}[thm]{Definition-Theorem}
\newtheorem{conjecture}[thm]{Conjecture}

\newcommand{\sF}{{\mathcal F}}

\newcommand{\sO}{{\mathcal O}}

\newcommand{\C}{{\mathbb C}}

\renewcommand{\P}{{\mathbb P}}
\newcommand{\Q}{{\mathbb Q}}
\newcommand{\R}{{\mathbb R}}

\newcommand{\ov}[1]{\overline{#1}}

\renewcommand{\S}{{\mathbb S}}

\newcommand{\Z}{{\mathbb Z}}

\newcommand{\qtq}[1]{\quad\mbox{#1}\quad}
\newcommand{\bp}{\bar{\partial}}
\newcommand{\Om}{\Omega}

\newcommand{\ts}{\otimes}

\newcommand{\btheorem}{\begin{theorem}}
\newcommand{\etheorem}{\end{theorem}}
\newcommand{\bproposition}{\begin{proposition}}
\newcommand{\eproposition}{\end{proposition}}
\newcommand{\bdefinition}{\begin{definition}}
\newcommand{\edefinition}{\end{definition}}
\newcommand{\bcorollary}{\begin{corollary}}
\newcommand{\ecorollary}{\end{corollary}}
\newcommand{\bproof}{\begin{proof}}
\newcommand{\eproof}{\end{proof}}
\newcommand{\bremark}{\begin{remark}}
\newcommand{\eremark}{\end{remark}}
\newcommand{\eexample}{\end{example}}
\newcommand{\bexample}{\begin{example}}

\newcommand{\elemma}{\end{lemma}}
\newcommand{\blemma}{\begin{lemma}}

\newcommand{\sq}{\sqrt{-1}}

\newcommand{\p}{\partial}

\renewcommand{\bar}{\overline}
\newcommand{\eps}{\varepsilon}

\renewcommand{\phi}{\varphi}

\newcommand{\ee}{\end{eqnarray*}}
\newcommand{\be}{\begin{eqnarray*}}

\newcommand{\beq}{\begin{equation}}
\newcommand{\eeq}{\end{equation}}

\newcommand{\bd}{\begin{enumerate}}
\newcommand{\ed}{\end{enumerate}}

\renewcommand{\tilde}{\widetilde}

\renewcommand{\>}{\rightarrow}

\usepackage{fancyhdr}
\renewcommand{\leftmark}{{\scriptsize Scalar curvature on compact complex manifolds}}

\begin{document}
\title{Scalar curvature on compact complex manifolds}
\makeatletter
\let\uppercasenonmath\@gobble
\let\MakeUppercase\relax
\let\scshape\relax
\makeatother
\author{Xiaokui Yang}
\date{}
\address{{Address of Xiaokui Yang: Morningside Center of Mathematics, Institute of
        Mathematics, Hua Loo-Keng Center of Mathematics,
        Academy of Mathematics and Systems Science,
        Chinese Academy of Sciences, Beijing, 100190, China.}}
\email{\href{mailto:xkyang@amss.ac.cn}{{xkyang@amss.ac.cn}}}
\maketitle

\maketitle

\begin{abstract}  In this paper, we prove that,  a compact
complex manifold $X$ admits a smooth Hermitian metric with positive
(resp. negative) scalar curvature if and only if $K_X$ (resp.
$K_X^{-1}$) is not pseudo-effective. On the contrary, we also show
that on an arbitrary compact
 complex
manifold $X$ with complex dimension $\geq 2$, there exist smooth
Hermitian metrics with positive \emph{total} scalar curvature, and
one of the key ingredients  in the proof relies on a recent solution
to the Gauduchon conjecture by G. Sz\'{e}kelyhidi, V. Tosatti and B.
Weinkove(\cite[Theorem~1.3]{STW}).

\end{abstract}

{\small{\setcounter{tocdepth}{1} \tableofcontents}}

\section{Introduction}

In this paper, we investigate the relationship between the sign of
the (total) scalar curvature of Hermitian metrics and the geometry
of the ambient complex manifolds.\\

On a compact K\"ahler manifold, one can define the positivity of
holomorphic bisectional curvature,  Ricci curvature, scalar
curvature and so on. The weakest one among them is the positivity of
total scalar curvature.  In algebraic geometry, the Kodaira
dimension can also characterize the positivity of the canonical
bundles and anti-canonical  bundles. In his seminal work
\cite{Yau2},
  Yau proved that,  on a  compact K\"ahler manifold $X$, if it admits a K\"ahler metric with  positive total scalar curvature,
  then the Kodaira dimension
 $\kappa(X)=-\infty$. Furthermore, Yau established that a compact
 K\"ahler surface is uniruled if and only if there exists a K\"ahler
 metric with positive total scalar curvature.    Recently, Heier-Wong pointed out in \cite{HW12} that a projective manifold is uniruled if it admits a K\"ahler metric with positive total scalar curvature.
 By using Boucksom-Demailly-Peternell-Paun's criterion for uniruled projective manifolds (\cite[Corollary~0.3]{BDPP13}),
 Chiose-Rasdeaconu-Suvaina obtained in \cite{CRS} a more general characterization that,  a compact Moishezon manifold is uniruled
if and only if it admits a smooth Gauduchon metric with positive
total Chern scalar curvature. As motivated by these
works(\cite{Yau2, BDPP13, HW12, CRS}), we investigate the total
Chern scalar curvature of Gauduchon metrics on general compact
complex manifolds.  Let $\omega$ be a smooth Hermitian metric on the
compact complex manifold $X$. For simplicity, we denote by $\mathscr
F(\omega)$ the total Chern scalar curvature of $\omega$, i.e.
$$\mathscr F(\omega)=\int_X s \omega^n=n
\int_X\text{Ric}(\omega)\wedge \omega^{n-1}.$$ Note that, when the
manifold is not K\"ahler, the total Chern scalar curvature differs
from the total scalar curvature of the Levi-Civita connection of the
underlying Riemannian metric (e.g. \cite{LY14}). Let $\mathscr W$ be
the space of smooth Gauduchon metrics on $X$. We obtain a complete
characterization on the image of the total scalar curvature function
$\mathscr F:\mathscr W\>\R$. (c.f. \cite[p.761-p.762]{HW12})

\btheorem\label{table} The image of the total scalar function
$\mathscr F:\mathscr W\>\R$ has exactly four different cases: \bd
\item $\mathscr F(\mathscr W)=\R$, if and only if neither $K_X$ nor
$K_X^{-1}$ is pseudo-effective;
\item $\mathscr F(\mathscr W)=\R^{>0}$, if and only if $K_X^{-1}$ is pseudo-effective but not unitary flat;
\item $\mathscr F(\mathscr W)=\R^{<0}$, if and only if $K_X$ is pseudo-effective but not unitary flat;
\item $\mathscr F(\mathscr W)=\{0\}$, if and only if $K_X$ is  unitary flat.
\ed \etheorem \noindent  One of the key ingredients in the proof of
Theorem \ref{table} relies on a border line case of Lamari's
positivity criterion (\cite{La99}) over compact complex manifolds
and Tosatti's characterizations for non-K\"ahler Calabi-Yau
manifolds (\cite{T}).
 Moreover, in Section
\ref{example} we exhibit a variety of \emph{non-K\"ahler Calabi-Yau}
manifolds which can distinguish all different cases in Theorem
\ref{table}.

\bremark More generally, for a Bott-Chern class $[\alpha]\in
H^{1,1}_{\mathrm{BC}}(X,\R)$, we can also define the function
$\mathscr F_{[\alpha]}:\mathscr W\>\R$ with respect to the class
$[\alpha]$
$$\mathscr F_{[\alpha]}(\omega)=\int_X [\alpha]\wedge \omega^{n-1},$$
which is well-defined since $\omega$ is Gauduchon. As analogous to
Theorem \ref{table}, we  obtain in Theorem \ref{class} a criterion
 for the positivity of the class $[\alpha]$.
It generalizes a result in \cite{Tel} which only deals with compact
complex surfaces.
 \eremark

 As an application of Theorem \ref{table} and Gauduchon's
conformal method, we obtain a criterion for the existence of
Hermitian metric with positive or negative (total) scalar curvature
on compact complex manifolds (see also \cite[Theorem~D]{CRS} for
some special cases):

\btheorem\label{scalar} Let $X$ be a compact complex manifold. The
following are equivalent \bd \item  $K_X$ (resp. $K_X^{-1}$) is not
pseudo-effective;

\item  $X$ has a
Hermitian metric with positive (resp. negative) scalar curvature;

\item $X$ has a
Gauduchon metric with positive (resp. negative) total scalar
curvature. \ed \etheorem

 On the other hand, it is well-known that if a compact complex
manifold $X$ admits a smooth Hermitian metric with positive scalar
curvature, then the Kodaira dimension $\kappa(X)=-\infty$. As
another application of Theorem \ref{table}, we obtain

\btheorem\label{main0} Let $X$ be a compact complex manifold. Then
there exists a smooth Gauduchon metric $\omega_G$ with vanishing
total scalar curvature if and only if $X$ lies in one of the
following cases \bd\item $\kappa(X)=-\infty$, and neither $K_X$ nor
$K^{-1}_X$ is pseudo-effective;
\item $\kappa(X)=-\infty$, and $K_X$ is unitary flat;
\item $\kappa(X)=0$, and  $K_X$ is a holomorphic torsion, i.e. $K_X^{\ts m}=\sO_X$ for some $m\in \Z^+$.\ed
 \etheorem

\noindent  It is easy to see that Theorem \ref{main0} excludes
non-K\"ahler Calabi-Yau manifolds with $\kappa(X)\geq 1$ (e.g.
Example \ref{kodaira}). More generally,
 one has (see also \cite[Proposition ~2.4]{ACS}):

\bcorollary\label{coro} Let $X$ be a compact complex manifold.
Suppose $\kappa(X)>0$, then for any Gauduchon metric $\omega$, the
total scalar curvature $\mathscr F(\omega)<0$. \ecorollary

\noindent The following result is a straightforward application of
Theorem \ref{main0}, and it appears to be new and interesting in its
own right.

 \bcorollary Let $X$ be a compact
K\"ahler manifold. If there exists a Gauduchon (e.g. K\"ahler)
metric $\omega$ with vanishing total scalar curvature, then either
\bd\item $\kappa(X)=-\infty$ and neither $K_X$ nor $K^{-1}_X$ is
pseudo-effective; or
\item $X$ is Calabi-Yau, i.e. there exists a smooth K\"ahler metric $\tilde \omega$ with $\emph{Ric}(\tilde \omega)=0$.\ed
\ecorollary

One may wonder whether similar results hold for Hermitian metrics
other than Gauduchon metrics. Unfortunately, one can not replace the
Gauduchon
 metric condition by an arbitrary Hermitian metric even if the ambient manifold is K\"ahler. More precisely, we
 obtain:

\btheorem\label{kahler}  Let $X$ be an arbitrary compact complex
manifold with $\dim X>1$. Then there exists a Hermitian metric
$\omega$ with positive total scalar curvature.  Moreover, if $X$ is
K\"ahler, then there exists a conformally K\"ahler metric $\omega$
with positive total scalar curvature (i.e. $\omega =e^{f}\omega_0$
for some K\"ahler metric $\omega_0$ and $f\in C^\infty(X,\R)$).

 \etheorem

\noindent The proof of Theorem \ref{kahler} relies on
Sz\'{e}kelyhidi-Tosatti-Weinkove's
 solution to the Gauduchon conjecture on compact complex manifolds (\cite[Theorem~1.3]{STW}, i.e. Theorem
 \ref{STW}).

\vskip 1\baselineskip

\noindent\textbf{Acknowledgement.} The author is very grateful to
Professor
 Kefeng Liu for his support, encouragement and stimulating
discussions over  years. He would like to thank Professor Valentino
Tosatti and the anonymous referee for their comments and suggestions
on an earlier versions of this paper which clarify and improve the
presentations. He would also like to thank Professors
 Shing-Tung Yau, Fangyang Zheng and Xiangyu Zhou for their interests and discussions.
 This work
was partially supported by China's Recruitment
 Program of Global Experts and National Center for Mathematics and Interdisciplinary Sciences,
 Chinese Academy of Sciences.

\vskip 2\baselineskip

\section{Preliminaries}

\subsection{Curvatures on complex manifolds}
 Let $(X,\omega_g)$ be a compact Hermitian manifold. The Chern connection
 on $(T^{1,0}X,\omega_g)$ has  Chern curvature components \beq R_{i\bar
j k\bar \ell}=-\frac{\p^2g_{k\bar \ell}}{\p z^i\p\bar z^j}+g^{p\bar
q}\frac{\p g_{k\bar q}}{\p z^i}\frac{\p g_{p\bar \ell}}{\p\bar
z^j}.\eeq  The (first Chern) Ricci form $\text{Ric}(\omega_g)$ of
$(X,\omega_g)$ has components
$$R_{i\bar j}=g^{k\bar \ell}R_{i\bar jk\bar \ell}=-\frac{\p^2\log\det(g)}{\p z^i\p\bar z^j}$$
and it is well-known that the Ricci form represents the first Chern
class of the complex manifold $X$. The (Chern) scalar curvature
$s_g$ of $(X,\omega_g)$ is defined as \beq
s_g=\text{tr}_{\omega_g}\text{Ric}(\omega_g)=g^{i\bar j} R_{i\bar
j}. \eeq The total scalar curvature is \beq \int_X s_g \omega_g^n=n
\int\text{Ric}(\omega_g)\wedge \omega_g^{n-1}.\eeq where $n$ is the
complex dimension of $X$. \bd
\item A Hermitian metric $\omega_g$ is called a Gauduchon metric if
$\p\bp\omega_g^{n-1}=0$. It is proved by Gauduchon (\cite{Ga2})
that, in the conformal class of each Hermitian metric, there exists
a unique Gauduchon metric (up to constant scaling).

\item A compact complex manifold $X$ is called a Calabi-Yau manifold if $c_1(X)=0\in
H^2(X,\R)$.

\item  A compact complex manifold $X$ is called uniruled if it is
covered by rational curves.\ed

\subsection{Positivity of line bundles.}

Let $(X,\omega_g)$ be a compact Hermitian manifold, and $L\>X$ be a
holomorphic line bundle. \bd\item $L$ is said to be \emph{positive}
(resp. \emph{semi-positive}) if there exists a smooth Hermitian
metric $h$ on $L$ such that the curvature form $R=-\sq\p\bp\log h$
is a positive (resp. semi-positive) $(1,1)$-form.

\item $L$ is said to be \emph{nef}, if for any  $\eps>0$, there exists a
smooth Hermitian metric $h_\eps$ on $L$ such that the curvature of
$(L,h_\eps)$ satisfies $ -\sq\p\bp\log h_\eps\geq -\eps \omega_g.$

\item $L$ is said to be \emph{pseudo-effective}, if there exists a
(possibly) singular Hermitian metric $h$ on $L$ such that the
curvature  of $(L,h)$ satisfies $-\sq\p\bp\log h\geq 0$ in the sense
of currents. (See \cite{Dem} for more details.)

\item $L$ is said to be \emph{ $\Q$-effective}, if there exists some
positive integer $m$ such that $H^0(X,L^{\ts m})\neq 0$.

\item $L$ is called \emph{unitary flat} if there exists a smooth Hermitian metric $h$ on $L$ such that the curvature of
$(L,h)$ is zero, i.e. $-\sq\p\bp\log h=0$.

 \item  The Kodaira dimension $\kappa(L)$ of $L$ is defined to
be
$$\kappa(L):=\limsup_{m\>+\infty} \frac{\log \dim_\C
H^0(X,L^{\ts m})}{\log m}$$ and the \emph{Kodaira dimension}
$\kappa(X)$ of $X$ is defined as $ \kappa(X):=\kappa(K_X)$ where the
logarithm of zero is defined to be $-\infty$. \ed

\subsection{Bott-Chern classes and Aeppli classes}
On  compact complex (especially, non-K\"ahler) manifolds, the
Bott-Chern
  cohomology and the Aeppli cohomology are very useful:
  $$ H^{p,q}_{\mathrm{BC}}(X):= \frac{\text{Ker} d \cap \Om^{p,q}(X)}{\text{Im} \p\bp \cap \Om^{p,q}(X)}\qtq{and} H^{p,q}_{\mathrm{A}}(X):=\frac{\text{Ker} \p\bp \cap \Om^{p,q}(X)}{\text{Im} \p \cap \Om^{p,q}(X)+ \text{Im}\bp \cap
  \Om^{p,q}(X)}.$$

\noindent Let $\mathrm{Pic}(X)$ be the set of holomorphic line
bundles over $X$. Similar to the first Chern class map
$c_1:\mathrm{Pic}(X)\>H^{1,1}_{\bp}(X)$, there is a \emph{first
Bott-Chern class} map \beq c_1^{\mathrm{BC}}:\mathrm{Pic}(X)\>
H^{1,1}_{\mathrm{BC}}(X).\eeq Given any holomorphic line bundle
$L\to X$ and any Hermitian metric $h$ on $L$, its curvature form
$\Theta_h$ is locally given by $-\sq\p\bp\log h$.  We define
$c_1^{\mathrm{BC}}(L)$ to be the class of $\Theta_h$ in
$H^{1,1}_{\mathrm{BC}}(X)$ (modulo a constant $2\pi$).  For a
complex manifold $X$, $c_1^{\mathrm{BC}}(X)$ is defined to be
$c_1^{\mathrm{BC}}(K^{-1}_X)$ where $K_X^{-1}$ is the anti-canonical
line bundle $\wedge^n T^{1,0}X$. It is easy to see that
$c_1^{\mathrm{BC}}(L)=0$ if and only if $L$ is unitary flat.

\vskip 2\baselineskip

\section{Total scalar curvature of Gauduchon metrics}

Let $X$ be a compact complex manifold of complex dimension $n$.
Suppose

\bd
\item[$\bullet$] $\mathscr E$ is the set of real $\p\bp$-closed $(n-1,n-1)$
forms on $X$;

\item[$\bullet$] $\mathscr V$ is the set of real positive $\p\bp$-closed $(n-1,n-1)$
forms on $X$;

\item[$\bullet$] $\mathscr G=\{\omega^{n-1}\ |\ \omega\ \text{is a Gauduchon metric
}\}$.

\ed

\noindent In \cite{Mi}, M.L. Michelsohn observed that the power map:
$\eta\>\eta^{n-1}$, from $\Lambda^{1,1}T_x^*X$ to
$\Lambda^{n-1,n-1}T_x^*X$, carries the cone of strictly positive
$(1, 1)$-forms bijectively onto the cone of strictly positive $(n-1,
n-1)$-forms at each point $x\in X$, and obtained

\blemma\label{michelson} $\mathscr V=\mathscr G$. \elemma

\noindent \emph{The proof of Theorem \ref{table}.}  \emph{Claim 1.}
The canonical bundle $K_X$ is pseudo-effective if and only if
$\mathscr F(\omega)\leq 0$ for every Gauduchon metric $\omega$. If
in addition, there exists some Gauduchon metric $\omega_0$ such that
$\mathscr F(\omega_0)=0$, then $\mathscr F(\mathscr W)=\{0\}$ and
$K_X^{-1}$ is unitary flat.

  Suppose $K_X$ is pseudo-effective, it is well-known that there exist a smooth Hermitian metric $\omega_1$ and a real valued
  function $\phi\in L^1(X,\R)$ such that
  $$\text{Ric}(\omega_1)+\sq\p\bp\phi\leq 0$$
in the sense of currents. Then for any smooth Gauduchon metric
$\omega\in\mathscr W$, \be \mathscr F(\omega)&=&n
\int_X\text{Ric}(\omega)\wedge \omega^{n-1}\\
&=&n\int_X\left(\text{Ric}(\omega_1)-\sq\p\bp\log\left(\frac{\omega^n}{\omega_1^n}\right)\right)\wedge\omega^{n-1}\\
&=&n\int_X\text{Ric}(\omega_1)\wedge\omega^{n-1}\\
&=&n\int_X\left(\text{Ric}(\omega_1)+\sq\p\bp\phi\right)\wedge
\omega^{n-1}\leq 0.\ee

Conversely, we assume $\mathscr F(\omega)\leq 0$ for every Gauduchon
metric $\omega$. Now we follow the strategy in \cite{La99} to show
$K_X$ is pseudo-effective. Suppose there exists some Gauduchon
metric $\omega_0$ such that $\mathscr F(\omega_0)=0$, we shall show
that $\mathscr F(\mathscr W)=\{0\}$ and there exists a Hermitian
metric $\tilde \omega$ with $\text{Ric}(\tilde \omega)=0$, i.e.
$K_X$ is unitary flat.

Indeed, for any fixed $\p\bp$-closed $(n-1,n-1)$ form $\eta\in
\mathscr E$, we define a real function $f:\R\>\R$ as
 $$ f(t)=n\int_X \text{Ric}(\omega_0)\wedge \left(
 (1-t)\omega_0^{n-1}+t\eta\right).$$ Since $\omega_0$ has zero total
 scalar curvature, we have
\beq f(t)=t\cdot \left(n\int_X \text{Ric}(\omega_0)\wedge
\eta\right).\label{linear}\eeq Hence, $f(t)$ is linear in $t$. On
the other hand, for small $|t|$, we have
$$(1-t)\omega^{n-1}_0+t\eta\in \mathscr V.$$
By Lemma \ref{michelson}, there exist a small number $\eps_0>0$ and
a family of Gauduchon metrics $\omega_t$ with
$t\in(-2\eps_0,2\eps_0)$ such that
$$\omega_t^{n-1}=(1-t)\omega^{n-1}_0+t\eta.$$
Then we have \be f(t)&=&n\int_X
\text{Ric}(\omega_0)\wedge\omega_t^{n-1}=n\int_X\left(\text{Ric}(\omega_t)-\sq\p\bp\log\left(\frac{\omega_0^n}{\omega_t^n}\right)\right)\wedge\omega_t^{n-1}\\
&=&n\int_X \text{Ric}(\omega_t)\wedge\omega_t^{n-1}=\mathscr
F(\omega_t)\leq 0,\ee for $t\in(-2\eps_0,2\eps_0)$. In particular,
we have  $f(\eps_0)\leq 0$ and $f(-\eps_0)\leq 0$. However, by
(\ref{linear}), $f(t)$ is linear in $t$ and $f(0)=0$. Hence,
$f(t)\equiv 0$, i.e.
$$n\int_X \text{Ric}(\omega_0)\wedge
\eta=0.$$ Since $\eta$ is an arbitrary element in $\mathscr E$, by
Lamari's  criterion \cite{La99}, there exists $\phi\in L^1(X,\R)$
such that $\text{Ric}(\omega_0)+\sq\p\bp \phi=0$ in the sense of
currents. Hence, $\phi\in C^\infty(X,\R)$ and the metric
$\tilde\omega=e^{-\frac{\phi}{n}}\omega_0 $ is Ricci-flat, i.e.
$$\text{Ric}(\tilde \omega)=-\sq\p\bp\log\tilde\omega^n=0.$$
 Therefore, $K_X$ is
unitary flat. For any Gauduchon metric $\omega$, we have
$$\mathscr F(\omega)=n\int_X
\text{Ric}(\omega)\wedge\omega^{n-1}=n\int_X\left(\text{Ric}(\tilde
\omega)-\sq\p\bp\log\left(\frac{\omega^n}{\tilde
\omega^n}\right)\right)\wedge\omega ^{n-1}=0.$$

 Next, we show if $\mathscr F(\omega)<0$ for every Gauduchon
metric $\omega$, then $K_X$ is pseudo-effective but not unitary
flat. It follows from Lemma \ref{michelson} and
\cite[Lemma~3.3]{La99}. Indeed, we  fix a smooth Gauduchon metric
$\omega_G$.  By Lemma \ref{michelson}, for any $\p\bp$-closed
positive $(n-1,n-1)$ form $\psi\in \mathscr V$, there exists a
smooth Gauduchon metric $\omega$ such that $\omega^{n-1}=\psi$. It
is easy to see that
$$\mathscr F(\omega)=n\int_X \text{Ric}(\omega)\wedge \omega^{n-1}=n\int_X \text{Ric}(\omega_G)\wedge \omega^{n-1}=n\int_X \text{Ric}(\omega_G)\wedge\psi<0 $$
By \cite[Lemma~3.3]{La99}, there exists $\phi\in L^1(X,\R)$ such
that
$$-\text{Ric}(\omega_G)+\sq\p\bp\phi\geq 0$$ in the sense of
currents. That means, $K_X$ is pseudo-effective.\\

\noindent \emph{Claim 2.} The anti-canonical bundle $K^{-1}_X$ is
pseudo-effective if and only if $\mathscr F(\omega)\geq 0$ for every
Gauduchon metric $\omega$. If in addition, there exists some
Gauduchon metric $\omega_0$ such that $\mathscr F(\omega_0)=0$, then
$\mathscr F(\mathscr W)=\{0\}$ and $K_X^{-1}$ is unitary flat.

 The proof of Claim $2$ is similar to that of Claim $1$.\\

\noindent \emph{Claim 3.} $\mathscr F(\mathscr W)=\R$ if and only if
neither $K_X$ nor $K_X^{-1}$ is pseudo-effective.

   Indeed, if there exist two Gauduchon metrics
$\omega_1$ and $\omega_2$ such that $\mathscr F(\omega_1)> 0$ and
$\mathscr F(\omega_2)<0$, then there exists a smooth Gauduchon
metric $\omega_G$ such that $\mathscr F(\omega_G)=0$, and by the
scaling relation
 \beq \mathscr F (\lambda \omega)=\lambda^{n-1}\mathscr F(\omega).\label{scaling}\eeq
we have $\mathscr F(\mathscr W)=\R$. Acually, by Lemma
\ref{michelson}, there exists a Gauduchon
 metric $\omega_G$ such that
 \beq \omega_G^{n-1}=\mathscr F(\omega_1)\omega_2^{n-1}-\mathscr F(\omega_2)\omega_1^{n-1}\eeq
and the total scalar curvature \be \mathscr F(\omega_G)&=&n\int_X
\text{Ric}(\omega_G)\wedge\omega_G^{n-1} =n\int_X
\text{Ric}(\omega_G)\wedge\left(\mathscr
F(\omega_1)\omega_2^{n-1}-\mathscr
F(\omega_2)\omega_1^{n-1}\right)\\
&=&\mathscr F(\omega_1)\cdot n\int_X
\text{Ric}(\omega_2)\wedge\omega_2^{n-1}-\mathscr F(\omega_2)\cdot
n\int_X
\text{Ric}(\omega_1)\wedge\omega_1^{n-1}\\
&=&0.\ee Now  Claim $3$ follows from Claim $1$ and Claim $2$. The
proof of Theorem \ref{table} is completed. \qed

\vskip 1\baselineskip

Before giving the proof of Theorem \ref{scalar}, we need the
following results which follow from Gauduchon's conformal method. We
refer to \cite{Balas1, Balas2, ACS, Ga2, Ga3} for more details. See
also an almost Hermitian version  investigated in \cite{LU}. For
readers' convenience, we include a proof here.

 \blemma\label{conformal} Let $X$ be a compact
complex manifold. The following are equivalent:

\bd \item  there exists a smooth Gauduchon metric with positive
(resp. negative,  zero) \emph{total} scalar curvature; \item  there
exists a smooth Hermitian metric with positive (resp. negative,
zero) scalar curvature.\ed \elemma

\bproof Let $\omega_G$ be a Gauduchon metric and $s_G$ be its Chern
scalar curvature. It is well-known (e.g. \cite{Ga2,Ga3}) that the
following equation \beq s_G- \text{tr}_{\omega_G}\sq\p\bp
f=\frac{\int_X s_G\omega_G^n}{\int_X\omega^n_G} \label{hopf}\eeq has
a solution $f\in C^\infty(X)$ since $\omega_G$ is Gauduchon and the
integration
$$\int_X\left(s_G-\frac{\int_X
s_G\omega_G^n}{\int_X\omega^n_G}\right)\omega_G^n=0.$$  Let $
\omega_g= e^{\frac{f}{n}}\omega_G$. Then the  (Chern) scalar
curvature $s_g$ of $\omega_g$ is, \be s_g&=&\text{tr}_{\omega_g
}\text{Ric}(\omega_g)=-\text{tr}_{\omega_g }\sq \p\bp\log(
\omega_g^n)\\&=&-f^{-\frac{1}{n}}\text{tr}_{ \omega_G}\sq\p\bp\log (e^f \omega_G^n) \\
&=&f^{-\frac{1}{n}}\left(s_G-\text{tr}_{\omega_G}\sq\p\bp
f\right)\\&=&f^{-\frac{1}{n}}\frac{\int_X
s_G\omega_G^n}{\int_X\omega^n_G}=\frac{f^{-\frac{1}{n}}}{\int_X\omega^n_G}\cdot
\mathscr F(\omega_G).\ee Hence, a smooth Gauduchon metric $\omega_G$
with positive (resp. negative, zero) total scalar curvature can
induce a smooth Hermitian metric with positive (resp. negative,
zero) scalar curvature.

 Conversely,  let
$\omega_G=f_0^{\frac{1}{n-1}}\omega$ be a Gauduchon metric in the
conformal class of $\omega$ for some strictly positive function
$f_0\in C^\infty(X)$. Let $s_G$ be the corresponding Chern scalar
curvature with respect to the Gauduchon metric $\omega_G$. Then we
obtain \begin{eqnarray} \int_Xs_G\omega_G^n\nonumber&=&n\int_X  \text{Ric}(\omega_G)\wedge \omega_G^{n-1}=n\int_X \text{Ric}(\omega)\wedge \omega_G^{n-1}\\
&=&n\int_X f_0 \text{Ric}(\omega)\wedge \omega^{n-1}=\int_X f_0
s\omega^n.
\end{eqnarray}
Hence, a Hermitian metric with positive (resp. negative, zero)
scalar curvature can induce a Gauduchon metric with positive (resp.
negative, zero) total scalar curvature
 \eproof

 By using standard Bochner technique (e.g.
\cite{T, Yang16}), it is easy to show that if a compact complex
manifold $X$ admits a smooth Hermitian metric with positive scalar
curvature, then the Kodaira dimension $\kappa(X)=-\infty$. Hence, by
Lemma \ref{conformal}, one has the well-known result

\bcorollary\label{vanishing} Let $X$ be a compact complex manifold.
Suppose $X$ has a smooth Gauduchon metric $\omega$ with positive
total scalar curvature, then $\kappa(X)=-\infty$. \ecorollary

\vskip 1\baselineskip

 \noindent \emph{The proof of Theorem
\ref{scalar}.} If  there exists a smooth Gauduchon metric with
positive total scalar curvature, by Theorem \ref{table}, we deduce
$K_X$ is not pseudo-effective. Conversely, if $K_X$ is not
pseudo-effective, then by Theorem \ref{table}, there exists a smooth
Gauduchon metric with positive total scalar curvature. Now Theorem
\ref{scalar} follows from Lemma \ref{conformal}.\qed

\vskip 1\baselineskip

 \noindent \emph{The proof of Theorem
\ref{main0}.} Suppose there exists some Gauduchon metric $\omega_0$
such that $\mathscr F(\omega_0)=0$, then by Theorem \ref{table},
$\mathscr F(\mathscr W)=\R$ or $\{0\}$. If  $\mathscr F(\mathscr
W)=\R$, then by Corollary, \ref{vanishing} $\kappa(X)=-\infty$ and
by Theorem \ref{table}, neither $K_X$ nor $K_X^{-1}$ is
pseudo-effective. On the other hand, if $\mathscr F(\mathscr
W)=\{0\}$, by Theorem \ref{table}, $K_X$ is unitary flat, i.e.
$c_1^{\mathrm{BC}}(X)=0$. If $\kappa(X)\geq 0$, then by
\cite[Theorem~1.4]{T}, $K_X$ is a holomorphic torsion, i.e.
$K_X^{\ts m}=\sO_X$ for some $m\in \Z^+$. Indeed, since $K_X^{-1}$
unitary flat and hence nef, suppose $0\neq \sigma \in H^0(X,K_X^{\ts
m})$, then $\sigma$ is nowhere vanishing
(\cite[Proposition~1.16]{DPS}), i.e. $K_X^{\ts m}=\sO_X$.  (Note
that there exist compact complex manifolds with $\kappa(X)=-\infty$
and $c_1^{\mathrm{BC}}(X)=0$, e.g. Example \ref{bad}).

Conversely, if $K_X$ is unitary flat, or neither $K_X$ nor
$K_X^{-1}$ is pseudo-effective, by Theorem \ref{table}, $0\in
\mathscr F(\mathscr W)$, i.e. there exists a smooth Gauduchon metric
with vanishing total scalar curvature.  \qed

\vskip 1\baselineskip

\noindent\emph{The proof of Corollary \ref{coro}.} Suppose $\mathscr
F(\mathscr W)=\R$ or $\R^{>0}$, by Corollary \ref{vanishing},
$\kappa(X)=-\infty$ which is a contradiction. Suppose $\mathscr
F(\mathscr W)=\{0\}$, then $c_1^{\mathrm{BC}}(X)=0$ and by Theorem
\ref{main0}, $\kappa(X)=-\infty$ or $\kappa(X)=0$ which is  a
contradiction again. Hence we have $\mathscr F(\mathscr
W)=\R^{<0}$.\qed

\vskip 1\baselineskip

\noindent By using similar ideas as in the proof of Theorem
\ref{table}, we obtain \btheorem\label{class} Let $X$ be a compact
complex manifold. For a Bott-Chern class $[\alpha]\in
H^{1,1}_{\mathrm{BC}}(X,\R)$, we define a function $\mathscr
F_{[\alpha]}:\mathscr W\>\R$ with respect to the class $[\alpha]$ as
\beq \mathscr F_{[\alpha]}(\omega)=\int_X [\alpha]\wedge
\omega^{n-1}.\eeq Then image of the total scalar function $\mathscr
F_{[\alpha]}:\mathscr W\>\R$ has four different cases: \bd
\item $\sF_{[\alpha]}(\mathscr W)=\R$, if and only if neither $[\alpha]$ nor
$-[\alpha]$ is pseudo-effective;
\item $\sF_{[\alpha]}(\mathscr W)=\R^{>0}$, if and only if $[\alpha]$ is pseudo-effective but not zero;
\item $\sF_{[\alpha]}(\mathscr W)=\R^{<0}$, if and only if $-[\alpha]$ is pseudo-effective but not zero;
\item $\sF_{[\alpha]}(\mathscr W)=\{0\}$, if and only if $[\alpha]$ is zero.
\ed \etheorem

\bremark Since $c_1^{\mathrm{BC}}(X)=c_1^{\mathrm{BC}}(K^{-1}_X)$,
if we set $[\alpha]= c_1^{\mathrm{BC}}(X)\in
H^{1,1}_{\mathrm{BC}}(X,\R)$ in Theorem \ref{class},  we establish
Theorem \ref{table}.

\eremark

\vskip 1\baselineskip

\section{Some open problems}\label{open}

Let $X$ be a compact complex manifold in class $\mathcal {C} $, i.e.
$X$ is bimeromorphic to a compact K\"ahler manifold. Compact
K\"ahler, Moishezon and projective manifolds are all in class
$\mathcal C$.

The following conjectures are either well-known or implicitly
indicated in the literatures in some special cases, and we refer to
 \cite{Yau2,BDPP13,HW12,CRS} and the references therein.
\begin{conjecture}\label{A} $\kappa(X)=-\infty$ if and only if $X$
is uniruled, i.e. $X$ is covered by rational curves.
\end{conjecture}

\begin{conjecture}\label{B} $K_X$ is pseudo-effective if and only if $K_X$ is
$\Q$-effective.
\end{conjecture}

\begin{conjecture}\label{C} $\kappa(X)=-\infty$ if and only if there exists a
Gauduchon metric with positive total scalar curvature.
\end{conjecture}

\bremark In Conjecture \ref{A},  Conjecture \ref{B} and Conjecture
\ref{C}, the necessary condition directions are well-known. \eremark

\bcorollary\label{eq} Conjecture \ref{B} is equivalent to Conjecture
\ref{C}. \ecorollary

\bproof Suppose Conjecture \ref{B} is valid. If $\kappa(X)=-\infty$
and for any Gauduchon metric $\omega$ the total scalar curvature is
non-positive, then by Theorem \ref{table} $K_X$ is pseudo-effective.
Hence by Conjecture \ref{B}, $K_X$ is $\Q$-effective, i.e.
$H^0(X,K_X^{\ts m})\neq 0$. That means $\kappa(X)\geq 0$ which is a
contradiction.

Assume Conjecture \ref{C} is true. Suppose $K_X$ is pseudo-effective
but $K_X$ is not $\Q$-effective, i.e. $\kappa(X)=-\infty$. By
Conjecture \ref{C},  there exists a Gauduchon metric with positive
total scalar curvature. According to Theorem \ref{table}, $K_X$ is
not pseudo-effective which is a contradiction. \eproof

\bcorollary If $X$ is Moishezon, Conjecture \ref{A},  Conjecture
\ref{B} and Conjecture \ref{C} are equivalent. \ecorollary

\bproof It follows by Corollary \ref{eq} and \cite[Theorem~D]{CRS}.
Indeed, it is shown in \cite[Theorem~D]{CRS} that,  a compact
Moishezon manifold $X$ is uniruled if and only if $X$ admits
Gauduchon metric with positive total scalar curvature. Hence,
Conjecture \ref{A} and Conjecture \ref{C} are equvialent.\eproof

%

The following conjecture is of particular interest in K\"ahler
geometry.
\begin{conjecture}\label{D} Let $X$ be a compact K\"ahler manifold.
Then $X$ is uniruled if and only if $X$ admits a smooth K\"ahler
metric with positive total scalar curvature.
\end{conjecture}

\noindent When $X$ is a compact K\"ahler surface, Conjecture \ref{D}
was proved by Yau (\cite[Section~1.2]{Yau2}). On the other hand,
Conjecture \ref{D} predicts that on compact K\"ahler manifolds, the
existence of rational curves requires merely the positivity of total
scalar curvature of some K\"ahler metric.

\vskip 1\baselineskip

\section{Examples of non-K\"ahler Calabi-Yau manifolds}\label{example}

In this section we give some examples of \emph{non-K\"ahler
Calabi-Yau} manifolds
 satisfying the conditions  in Theorem \ref{table} or Theorem
\ref{main0}. These examples also show  significant differences
between non-K\"ahler manifolds and K\"ahler manifolds in our
setting.
\subsection{}
Let $X=\S^{2n-1}\times \S^1$ be the standard $n$-dimensional ($n\geq
2$) Hopf manifold. It is diffeomorphic to $\C^n- \{0\}/G$ where $G$
is cyclic group generated by the transformation $z\rightarrow
\frac{1}{2}z$.   On $X$, there is a natural induced Hermitian metric
$\omega$ given by
 $$ \omega=\sq h_{i\bar j}dz^i\wedge d\bar z^j=
\frac{4\delta_{i\bar j}}{|z|^2}dz^i\wedge  d\bar z^j.$$ This example
is studied with details in \cite{LY12, LY14, T, TW13, TW, Yang17A,
Yang17B}. One has
$$ \mathrm{Ric}(\omega)=-\sq\p\bp\log\det\omega^n=n\cdot\sq
\p\bp\log |z|^2.$$ Hence $\mathrm{Ric}(\omega)$ is semi-positive. In
particular, $K^{-1}_X$ is pseudo-effective. Moreover, for any
Gauduchon metric $\omega_G$, we have $\mathscr F(\omega_G)>0$.
Indeed, \be \mathscr F(\omega_G)&=&n\int_X
\mathrm{Ric}(\omega_G)\wedge\omega_G^{n-1}\\
&=&n\int_X\left(
\mathrm{Ric}(\omega)-\sq\p\bp\log\left(\frac{\omega_G^n}{\omega^n}\right)\right)\wedge\omega_G^{n-1}\\
&=&n\int_X \mathrm{Ric}(\omega)\wedge\omega_G^{n-1}=\int_X
\left(\mathrm{tr}_{\omega_G}\mathrm{Ric}(\omega)\right)\omega_G^{n}>0.\ee
Moreover, $X$ contains no rational curve. Otherwise, we have a
nonzero holomorphic map from $\P^1$ to $\C^n$ which is absurd. It is
easy to see $K_X$ is not pseudo-effective. In summary,

\bexample\label{Hopfsurface} On Hopf manifold $X=\S^{2n-1}\times
\S^1$,  we have the following properties:

\bd \item $\kappa(X)=-\infty$;

\item the total scalar curvature
$\mathscr F(\omega)>0$ for every Gauduchon metric $\omega$;

\item $X$ has no rational curve; (Counter-example to Conjecture \ref{A} on general complex manifolds)
\item $c_1(X)=0$ and $c_1^{\mathrm {BC}}(X)\neq 0$;
\item $K^{-1}_X$ is pseudo-effective, and $K_X$ is not pseudo-effective.
\ed This example lies in  case $(2)$ of Theorem \ref{table}.
\eexample

\subsection{} We give an example described  in \cite{T, Ma}. Let
$\alpha,\beta$ be the two roots of the equation $x^2-(1+i)x+1=0.$
The minimal polynomial over $\Q$ of $\alpha$ (and $\ov{\beta}$) is
$x^4-2x^3+4x^2-2x+1.$ Let $\Lambda$ be the lattice in $\C^2$ spanned
by the vectors $(\alpha^j,\ov{\beta}^j)$, $j=0,\dots,3$. Let
$Y=\C^2/\Lambda$. The automorphism of $\C^2$ given by multiplication
by $\begin{pmatrix}
\alpha & 0  \\
0 & \ov{\beta}
\end{pmatrix}$
descends to an automorphism $f$ of $Y$.  Let
$C=\mathbb{C}/(\mathbb{Z}\oplus\mathbb{Z}\tau)$ be an elliptic
curve, and we define a holomorphic free $\mathbb{Z}^2$-action on
$Y\times\C$ by
$$(1,0)\cdot(x,z)=(x,z+1),\quad (0,1)\cdot(x,z)=(f(x),z+\tau).$$
Then the quotient space  $X$ is a holomorphic fiber bundle $X\to C$
with fiber $Y$. Following \cite{T, Ma}, $X$ is a non-K\"ahler
manifold with $c_1^{\mathrm{BC}}(X)=0$ and $\kappa(X)=-\infty$.

\bexample\label{bad} On $X$,   we have the following properties:

 \bd
 \item $\kappa(X)=-\infty$;

 \item $c_1(X)=0$, $c^{\mathrm{BC}}_1(X)=0$, i.e. $K_X$ is unitary flat;
\item the total scalar curvature $\mathscr F(\omega)\equiv0$
for every Gauduchon metric $\omega$;

\item $X$ has no rational curve;

\item $K_X$ is pseudo-effective, but it is not $\Q$-effective. (Counter-example to Conjecture \ref{B} on general complex
manifolds). \ed This example lies in  case $(4)$ of Theorem
\ref{table} and case $(2)$ of Theorem \ref{main0}. \eexample

\subsection{} This construction follows from \cite[Example~3.4]{T}. Let $T=\mathbb{C}^n/\Lambda$ be an
$n$-torus, $\Sigma$ be a compact Riemann surface of genus $g\geq 2$
and $\pi:X\to\Sigma$ be any topologically nontrivial principal
$T$-bundle over $\Sigma$. It is shown in \cite[Example~3.4]{T} that
$X$ is a non-K\"ahler manifold with $c_1(X)=0$,
$c^{\mathrm{BC}}_1(X)\neq 0$ and $\kappa(X)=1$. By Corollary
\ref{coro}, we know for any Gauduchon metric $\omega$, the  total
scalar curvature $\mathscr
 F(\omega)<0$.
 \bexample\label{kodaira}  On $X$,   we
have the following properties:

 \bd
 \item $\kappa(X)=1$;

\item $c_1(X)$=0 and $c^{\mathrm{BC}}_1(X)\neq 0$;

\item the total scalar curvature $\mathscr F(\omega)<0$
for every Gauduchon metric $\omega$;

\item $K_X$ is pseudo-effective, but $-K_X$ is not pseudo-effective.

\ed This example lies in  case $(3)$ of Theorem \ref{table}.
\eexample

\bexample\label{4} Let $X=X_2\times X_3$ be the product manifold
where $X_2$ and $X_3$ are the complex manifolds constructed in
Example \ref{bad} and Example \ref{kodaira} respectively. It has the
following properties:

 \bd
 \item $\kappa(X)=-\infty$;

 \item $c_1(X)=0$, $c^{\mathrm{BC}}_1(X)\neq 0$.
\item the total scalar curvature $\mathscr F(\omega)<0$
for every Gauduchon metric $\omega$;
%
%
%
%

\item $K_X$ is pseudo-effective, but it is not $\Q$-effective.  \ed
This example lies in  case $(3)$ of Theorem \ref{table}. \eexample

 \bexample Let $X=X_1\times X_2 \times X_3$ be the product
manifold where $X_1$,  $X_2$ and $X_3$ are the complex manifolds
constructed in Example \ref{Hopfsurface}, Example \ref{bad} and
Example \ref{kodaira} respectively. It has the following properties:

 \bd
 \item $\kappa(X)=-\infty$.

 \item $c_1(X)=0$, $c^{\mathrm{BC}}_1(X)\neq 0$.
\item The total scalar curvature can be any real number. Indeed, it
follows from Example \ref{Hopfsurface} and Example \ref{4} by using
the scaling trick.
%
%
%
%

\item  Neither $K_X$ nor $K_X^{-1}$ is pseudo-effective.  \ed
This example lies in  case $(1)$ of Theorem \ref{table} and case
$(1)$ of Theorem \ref{main0}.

\eexample

\bexample Let $X$ be a Kodaira surface (a non-K\"ahler compact
complex surface with torsion canonical line bundle).  It has the
following properties:

 \bd
 \item $\kappa(X)=0$;

 \item  $c_1(X)=0$, $c^{\mathrm{BC}}_1(X)=0$ and $K_X$  is unitary
 flat;
\item  the total scalar curvature $\mathscr F(\omega)=0$
for every Gauduchon metric $\omega$;
 \ed
This example lies in  case $(4)$ of Theorem \ref{table} and case
$(3)$ of Theorem \ref{main0}.

\eexample

\vskip 2\baselineskip

\section{Existence of smooth Hermitian metrics with positive total scalar curvature}

In this section we prove Theorem \ref{kahler}. As we pointed out
before, one of the key gradients in the proof is
Sz\'{e}kelyhidi-Tosatti-Weinkove's
 solution to the Gauduchon conjecture  on compact complex manifold, which
is analogous to Yau's solution to the Calabi conjecture \cite{Yau78}
on compact  K\"ahler manifolds:

\btheorem\cite[Theorem~1.3]{STW}\label{STW} Let $X$ be a compact
complex manifold. Let $\omega_0$ be a smooth  Gauduchon metric, and
$\Phi$ be a closed real $(1, 1)$ form on $X$ with $[\Phi] =
c_1^{\mathrm{BC}} (X) \in H^{1,1}_{\mathrm{BC}}(X, \R)$. Then there
exists a smooth Gauduchon metric $\omega$ satisfying $[\omega^{ n-1}
] = [\omega_0^{n-1}]$ in $H^{n-1,n-1}_{\mathrm{A}} (X, \R)$ and \beq
\emph{Ric}(\omega) = \Phi.\eeq In particular, for any smooth volume
form $\sigma$ on $X$, there exists a smooth Gauduchon metric
$\omega$ such that \beq \omega^n=\sigma.\eeq

\etheorem

\btheorem  Let $X$ be an arbitrary compact complex manifold with
$\dim X>1$. Then there exists a smooth Hermitian metric $\omega$
with positive total scalar curvature.  Moreover, if $X$ is K\"ahler,
then there exists a conformally K\"ahler metric $\omega$ with
positive total scalar curvature (i.e. $\omega =e^{f}\omega_0$ for
some K\"ahler metric $\omega_0$ and $f\in C^\infty(X,\R)$).

\etheorem

\bproof There exists a smooth Gauduchon metric $\omega$ such that
the scalar curvature of $\omega$ is strictly positive at some point
$p\in X$. Indeed, fix an arbitrary smooth Hermitian metric
$\omega_0$ on $X$, then there exists some smooth function $F$ such
that
$$\text{Ric}(\omega_0)-\sq\p\bp F$$
is positive definite at point $p$. Let $\{z^i\}$ be the holomorphic
coordinates centered at point $p$, we can  choose $F(p)=-\lambda
|z|^2+o(|z|^3)$ for some large positive constant $\lambda$
(depending on $\text{Ric}(\omega_0)(p)$). On the other hand, by
Theorem \ref{STW}, there exists a smooth Gauduchon metric $\omega$
such that \beq \omega^n =e^F\omega_0^n.\label{CY}\eeq Then
$\text{Ric}(\omega)=\text{Ric}(\omega_0)-\sq\p\bp F$ is positive
definite at point $p$. Hence the scalar curvature of $\omega$ is
positive at  $p$. Let $\omega_f=e^{f}\omega$, then the scalar
curvature $s_f$ of $\omega_f$ is
$$s_f=\text{tr}_{\omega_f}\text{Ric}(\omega_f)=e^{-f}\text{tr}_{\omega}\left(\text{Ric}(\omega)-n\sq\p\bp f\right)=e^{-f}(s_\omega-n\Delta_\omega f)$$
and the total scalar curvature of $\omega_f$ is \be \int_X s_f
\omega_f^n&=&\int_X e^{(n-1)f}\cdot
(s_\omega-n\Delta_\omega f)\cdot \omega^n\\
&=&\int_X e^{(n-1)f}\cdot s_\omega\cdot
\omega^n-n^2\int_Xe^{(n-1)f}\cdot\sq\p\bp f\wedge \omega^{n-1}\\
&=&\int_X e^{(n-1)f}\cdot s_\omega\cdot
\omega^n+n^2(n-1)\int_Xe^{(n-1)f}\cdot\sq\p f\wedge \bp f\wedge
\omega^{n-1}\\
&&-n^2\int_Xe^{(n-1)f}\sq\bp f\wedge \p\omega^{n-1}, \ee where we
use Stokes' theorem in the last identity. It is obvious that
$$n^2(n-1)\int_Xe^{(n-1)f}\cdot\sq\p f\wedge \bp f\wedge
\omega^{n-1}\geq 0.$$ Moreover, we have \be
-n^2\int_Xe^{(n-1)f}\sq\bp f\wedge
\p\omega^{n-1}&=&-\frac{n^2}{n-1}\int_X\sq \left(\bp
e^{(n-1)f}\right) \wedge \p\omega^{n-1}\\
&=&\frac{n^2}{n-1}\int_X\sq  e^{(n-1)f} \cdot \bp\p\omega^{n-1}\\
&=&0\ee since $\omega$ is Gauduchon. Therefore, we obtain \beq
\int_X s_f \omega_f^n \geq \int_X e^{(n-1)f}\cdot s_\omega\cdot
\omega^n. \eeq
 Since
$s_\omega(p)>0$, by standard analytic techniques, there exists some
smooth function $f$ such that
$$\int_X s_f \omega_f^n\geq \int_X e^{(n-1)f}\cdot s_\omega\cdot \omega^n>0.$$
Indeed, without loss of generality, we can assume \beq \max_X
s_\omega=k_0>0,\ \ \ \min_X s_\omega=-k_1<0 \qtq{and} \int_X
\omega^n=1.\eeq Let
$$X_1=\left\{q\in X \ |\ s_\omega(q)\in\left[\frac{k_0}{2}, k_0\right]\right\},\ \ X_2=\left\{q\in X \ |\ s_\omega(q)\in\left[\frac{k_0}{4}, \frac{k_0}{2}\right)\right\}$$
and
$$X_3=\left\{q\in X\ |\ s_\omega(q)\in\left[-k_1, \frac{k_0}{4}\right)\right\}.$$
Let $f$ be a smooth function such that $e^{(n-1)f}\equiv 1$ on $X_3$
and  $e^{(n-1)f}\equiv 1+\frac{2k_1}{k_0}$ on $X_1$. Then \be \int_X
e^{(n-1)f}\cdot s_\omega\cdot \omega^n&\geq&  \int_{X_1}
e^{(n-1)f}\cdot s_\omega\cdot \omega^n+\int_{X_3} e^{(n-1)f}\cdot
s_\omega\cdot \omega^n\\
&\geq &-k_1+\left(1+\frac{2k_1}{k_0}\right)\cdot
\frac{k_0}{2}=\frac{k_0}{2}>0.\ee The proof of the first part of
Theorem \ref{kahler} is completed.

 Suppose $X$ is K\"ahler, then by the same arguments as above, there exists a smooth K\"ahler metric $\omega$
such that the scalar curvature of $\omega$ is strictly positive at
some point $p\in X$ where we use the Calabi-Yau theorem in
(\ref{CY}). By using the conformal method and integration by parts,
we obtain a conformally K\"ahler metric with positive total scalar
curvature. \eproof


\end{document}